\documentclass{article}
\usepackage{enumerate}
\usepackage{amsmath}
\usepackage{amssymb}
\usepackage{amsthm}
\DeclareSymbolFontAlphabet{\Bbb}{AMSb}
\newcommand{\id}{\text{{\rm id}}}
\newcommand{\Span}{\text{{\rm span}}}

\newcommand{\norm}[1]{ \|#1 \| }

\newcommand{\cotype}[1]{\mathbf{C_{#1}}}

\newcommand{\N}{\Bbb{N}}

\newcommand{\LL}{{\cal L}}

\newcommand{\Id}{\hookrightarrow}

\newcommand{\ui}{\mathcal{S}}
\theoremstyle{definition}
\newtheorem{definition}{Definition}

\theoremstyle{plain}

\newtheorem{theorem}[definition]{Theorem}
\newtheorem{corollary}[definition]{Corollary}
\newtheorem{proposition}[definition]{Proposition}
\theoremstyle{remark}
\numberwithin{equation}{section}
\begin{document}
\bibliographystyle{amsalpha}
\title{\bf $\boldsymbol{\Lambda(p)}$-sets and the limit order of operator 
ideals}
\author{Carsten Michels \\
\small e-mail: {\tt michels}@{\tt mathematik.uni-oldenburg.de}}
\date{}
\maketitle
\begin{abstract}
Given an infinite set $\Lambda$ of characters on a compact abelian group
 we show that $\Lambda$ is a $\Lambda(p)$-set for all $2<p<\infty$ if and 
only if the limit order of the ideal $\Pi_\Lambda$ of all $\Lambda$-summing
 operators coincides with that of the ideal $\Pi_\gamma$ of all 
Gaussian-summing operators, i.\,e. $\lambda(\Pi_\Lambda,u,v)=\lambda(
\Pi_\gamma,u,v)$ for all $1 \le u,v \le \infty$. This is a natural 
counterpart to a recent result of Baur which says that $\Lambda$ is a 
Sidon set if and only if $\Pi_\Lambda=\Pi_\gamma$.
Furthermore, our techniques,
 which are mainly based on complex interpolation, lead us to exact asymptotic
 estimates of the Gaussian-summing norm $\pi_\gamma(\id: \ui_u^n \Id 
\ui_v^n)$ of identities between finite-dimensional Schatten classes 
$\ui_u^n$ and $\ui_v^n$, $1 \le u,v \le \infty$. 
\end{abstract}
\section{Introduction and results}
We use standard notation and notions from Banach space theory, as presented
 e.\,g. in \cite{djt}, \cite{lt} and \cite{tj}. If $E$ is a Banach space, 
then $B_E$ is its (closed) unit ball and $E'$ its dual; we consider 
complex Banach spaces only. As usual $\LL(E,F)$ denotes the Banach space
 of all (bounded and linear) operators from $E$ into $F$ endowed with the 
operator norm.
\par
For an infinite orthonormal system $B \subset L_2(\mu)$ (over some 
probability space $(\Omega,\mu)$) an operator $T:E \rightarrow F$ between
 Banach spaces $E$ and $F$ is said to be $B$-summing if there exists a 
constant $c>0$ such that for all finite sequences
 $b_1, \ldots, b_n$ in $B$ and $x_1, \ldots,x_n$ in $E$
\begin{equation}
\label{bpsumming}
\left( \int_\Omega \norm{\sum_{i=1}^n b_i \cdot Tx_i}^2 d\mu \right)^{1/2}
 \le c  \cdot \sup_{x' \in B_{E'}} \left ( \sum_{i=1}^n |\langle 
x',x_i \rangle|^2 \right )^{1/2};
\end{equation}
we write $\pi_B(T)$ for the smallest constant $c$ satisfying
\eqref{bpsumming}. In this way we obtain the injective and maximal Banach operator ideal 
$(\Pi_B, \pi_B)$, which became of interest recently in the theses of 
Baur~\cite{baurdiss} and Seigner~\cite{seigner}.
For a 
sequence of independent standard Gaussian random variables 
 the associated Banach operator ideal $\Pi_\gamma$ of all Gaussian-summing 
operators was introduced by Linde and Pietsch \cite{linde}, and for 
 operators acting on finite-dimensional Hilbert spaces 
$\pi_\gamma$ is also known as the $\ell$-norm, which turned out to be 
important for the study of the geometry of Banach spaces (see e.\,g. 
\cite{tj}). 
\par
For an infinite subset $\Lambda$ of the character group $\Gamma$ of 
some compact abelian group $G$ (which can be viewed as an orthonormal 
system in $L_2(G,m_G)$, where $m_G$ denotes the normalized Haar measure 
on $G$) Baur in \cite[9.5]{baurdiss} (see also \cite[4.2]{baur99}) gave
 the following characterization:
\begin{center}
\em 
$\Lambda$ is a Sidon set if and only if $\Pi_\Lambda=\Pi_\gamma$.
\end{center}
Recall that a subset $ \Lambda \subset
\Gamma$ is said to be a {\em Sidon set}
 if there exists $\theta>0$ such that for all $Q=\sum_{\gamma \in \Lambda}
 \alpha_\gamma \cdot \gamma \in \Span(\Lambda)$ we have 
$
\sum_{\gamma \in \Lambda} |\alpha_\gamma| \le \theta \cdot 
\norm{Q}_{L_\infty(G,m_G)},
$
and for $2 < p < \infty$ it is called a {\em $\Lambda(p)$-set} if 
there exists a constant $c>0$ such that for all 
$\lambda \in \Span (\Lambda)$ we have $\norm{\lambda}_{L_p(G,m_G)} \le 
\norm{\lambda}_{L_2(G,m_G)}$; the infimum over all such constants $c$ 
is denoted by $K_p(\Lambda)$. Pisier in \cite{pisier78} showed that
 $\Lambda$ is a Sidon set if and only if
\begin{center}
\em 
 $\Lambda$ is a $\Lambda(p)$-set
 with $K_p(\Lambda) \le \kappa \sqrt{p}$ for all $2<p<\infty$ and some 
$\kappa>0$.
\end{center}
 As a natural counterpart of Baur's result we prove the 
following characterization of sets which are $\Lambda(p)$-sets for all 
$2<p<\infty$, with no control of $K_p(\Lambda)$ as in Pisier's 
characterization of Sidon sets above:
\begin{theorem}
\label{khinchar}
For every infinite subset $\Lambda \subset \Gamma$ the following are 
equivalent:
\begin{enumerate}[(a)]
\item
$\Lambda$ is a $\Lambda(p)$-set for all $2 < p < \infty$.
\item
$ \lambda(\Pi_\Lambda,u,v) = \lambda(\Pi_\gamma,u,v)$ for all $1 \le u,v \le
 \infty$.
\end{enumerate}
\end{theorem}
Here, the limit order $\lambda({\mathcal{A}},u,v)$ of a Banach operator ideal
 $({\mathcal{A}},A)$ for $1 \le u,v \le \infty$ is defined as usual 
(see e.\,g. \cite[14.4]{pietsch}): 
$$
\lambda({\mathcal{A}},u,v):= 
\inf \{ \lambda>0 \,| \, \exists \, \rho >0 \, \forall \, n \in \N:
 A(\id: \ell_u^n \Id \ell_v^n) \le \rho \cdot n^\lambda \}.
$$
Note that there exist sets which are $\Lambda(p)$-sets for all $2<p<\infty$
 but fail to be Sidon sets (see e.\,g. \cite[5.14]{lopez}). Our proof is
 mainly based on complex interpolation techniques, in particular on 
formulas for the complex interpolation of spaces of operators due to 
Pisier~\cite{pisier90} and Kouba~\cite{kouba} (see also \cite{dm99}). 
These techniques also yield asymptotic estimates of the Gaussian-summing
 norm of identities between finite-dimensional Schatten classes:
\begin{theorem}
\label{gammaschatten}
For $1 \le u,v \le \infty$
$$
\pi_\gamma(\ui_u^n \Id \ui_v^n) \asymp
\begin{cases}
n^{1/2+1/v} &\text{ if } 2 \le u \le \infty, \\
n^{1/2+ \max(0,1/2+1/v-1/u)} &\text{ if } 1 \le u \le 2. 
\end{cases}
$$
In particular,
$$
\pi_\gamma(\ui_u^n \Id \ui_v^n) \asymp n^{1/2 + \lambda(\Pi_\gamma,u,v)}.
$$
\end{theorem}
Here, $\ui_u^n$ denotes the space of all linear operators $T:\ell_2^n 
\rightarrow \ell_2^n$ endowed with the norm $\norm{T}_{\ui_u^n}:=
\norm{(s_k(T))_{k=1}^n}_{\ell_u^n}$, where $(s_k(T))_{k=1}^n$ is the sequence
 of singular numbers of $T$. Besides the interpolation argument, 
our proof uses the close relationship between the 
Gaussian-summing norm of the identity operator $\id_E$
 and the Dvoretzky dimension of a finite-dimensional Banach space $E$ 
due to Pisier.
\section{Complex interpolation of $\boldsymbol{B}$-summing operators}
Our main tool will be an ``interpolation theorem'' for the 
$B$-summing norm of a fixed operator acting between  
finite-dimensional complex interpolation
 spaces; a similar approach for the $(s,2)$-summing norm was used in
 \cite{dm98} to study the well-known ``Bennett--Carl inequalities'' 
within the context of interpolation theory.
\par 
For all information on complex interpolation we refer to \cite{BL}. 
Given an interpolation couple $[E_0,E_1]$ of complex Banach spaces and 
$0<\theta<1$, the associated complex interpolation space is denoted by 
$[E_0,E_1]_\theta$. If $E_0$ and $E_1$ are finite-dimensional Banach spaces
 with the same dimensions, we speak of a finite-dimensional interpolation 
couple, and in this case we define
$$
d_\theta[E_0,E_1]:=\sup_m \norm{\LL(\ell_2^m,[E_0,E_1]_\theta) \Id 
[\LL(\ell_2^m,E_0),\LL(\ell_2^m,E_1)]_\theta}.
$$
Pisier~\cite{pisier90} and Kouba~\cite{kouba} derived upper estimates 
for $d_\theta[E_0,E_1]$ for particular situations (see also \cite{dm99});
 we will use the fact that for $1 \le p_0,p_1 \le 2$ 
\begin{equation}
\label{dtheta1}
d_\theta[\ell_{p_0}^n,\ell_{p_1}^n] \le \sqrt{2};
\end{equation}
in particular, $\sup_n d_\theta[\ell_1^n,\ell_2^n] < \infty$.
Junge \cite[4.2.6]{junge} gave an analogue for Schatten classes:
\begin{equation}
\label{dtheta2}
\sup_n d_\theta[\ui_1^n,\ui_2^n] <\infty.
\end{equation}
These ``uniform'' estimates will be crucial for the applications of the 
following result:
\begin{proposition}
\label{interpolB}
Let $0 < \theta <1$. Then for two 
finite-dimensional interpolation couples $[E_0,E_1]$ and $[F_0,F_1]$, each 
$T \in \LL([E_0,E_1]_\theta,[F_0,F_1]_\theta)$ and each orthonormal system 
$B \subset L_2(\mu)$ 
$$
\pi_B(T:[E_0,E_1]_\theta  \rightarrow [F_0,F_1]_\theta) 
\le d_\theta[E_0,E_1] \cdot  \pi_B(T:E_0 \rightarrow F_0)^{1-\theta}
 \cdot \pi_B(T:E_1 \rightarrow F_1)^\theta.
$$
\end{proposition}
\proof For the moment set $E_\theta:=[E_0,E_1]_\theta$, $F_\theta:=
[F_0,F_1]_\theta$, and
consider for $\eta=0,\theta,1$ and ${\mathcal{F}}= \{ b_1, \ldots, 
b_m \} \subset B$  the mapping
$$
\begin{array}{lccc}
\Phi_\eta^{m,{\mathcal{F}}} : 
& \LL(\ell_2^m,E_\eta) & \rightarrow & L_2(\mu,F_\eta)
 \\ & S  & \mapsto  &  \sum_{i =1}^m b_i \cdot TSe_i.
\end{array}
$$
Since for each $S= \sum_{j=1}^m e_j \otimes x_j \in \LL(\ell_2^m, E_\eta)$
$$
\norm{S}= \sup_{x' \in B_{E_\eta'}} \left ( \sum_{j=1}^m |\langle x',x_j 
\rangle|^2 \right )^{1/2},
$$
we obviously get that 
$$
\pi_B(T: E_\eta \rightarrow F_\eta) = \sup \{\norm{\Phi_\eta^{m,
 {\mathcal{F}}}} \, | \, m \in \N,  {\mathcal{F}} 
\subset B \text{ with }  |{\mathcal{F}}| = m \}.
$$ 
For the interpolated mapping 
$$
[\Phi_0^{m,{\mathcal{F}}}, \Phi_1^{m,{\mathcal{F}}}]_\theta : 
[\LL(\ell_2^m,E_0), \LL(\ell_2^m,E_1)]_\theta \rightarrow
[L_2(\mu,F_0), L_2(\mu,F_1)]_\theta
$$ 
by the usual interpolation theorem 
$$
\norm{[\Phi_0^{m,{\mathcal{F}}}, \Phi_1^{m,{\mathcal{F}}}]_\theta} \le 
\norm{\Phi_0^{m,{\mathcal{F}}}}^{1-\theta} \cdot 
\norm{\Phi_0^{m,{\mathcal{F}}}}^\theta.
$$
Since $[L_2(\mu,F_0),L_2(\mu,F_1)]_\theta = L_2(\mu,[F_0,F_1]_\theta)$ 
(isometrically, see \cite[5.1.2]{BL}) we obtain
$$
\norm{\Phi_\theta^{m,{\mathcal{F}}}} \le \norm{\LL(\ell_2^m,
[E_0,E_1]_\theta)
 \Id [\LL(\ell_2^m,E_0),\LL(\ell_2^m,E_1)]_\theta } \cdot 
\norm{[\Phi_0^{m,{\mathcal{F}}}, \Phi_1^{m,{\mathcal{F}}}]_\theta} .
$$
Consequently
\begin{align*}
\pi_B(T & :  [E_0,E_1]_\theta \rightarrow [F_0,F_1]_\theta) \\ 
&= \sup \{\norm{\Phi_\theta^{m, {\mathcal{F}}}} \, | \, m \in \N,  {\mathcal{F}} 
\subset B \text{ with }  |{\mathcal{F}}| = m \} \\
& \le  \sup 
\{d_\theta[E_0,E_1] \cdot \norm{\Phi_0^{m,
 {\mathcal{F}}}}^{1-\theta} \cdot \norm{\Phi_0^{m, {\mathcal{F}}}}^\theta
 \, | \, m \in \N,  {\mathcal{F}} 
\subset B \text{ with }  |{\mathcal{F}}| = m \} \\
& \le d_\theta[E_0,E_1] \cdot 
\pi_B(T:E_0 \rightarrow F_0)^{1-\theta} \cdot 
\pi_B(T: E_1 \rightarrow F_1)^\theta,
\end{align*}
the desired result. 
\qed
\par
Since for $0<\theta<1$ and $1 \le p_0,p_1,p_\theta \le \infty$ such that
 $1/{p_\theta} = (1-\theta)/{p_0} + \theta/{p_1}$ it holds 
$[\ell_{p_0}^n,\ell_{p_1}^n]_\theta= \ell_{p_\theta}^n$ isometrically 
(see \cite[5.1.1]{BL}), we obtain together with \eqref{dtheta1} the 
 following corollary:
\begin{corollary}
\label{interpollimit}
For $0<\theta<1$ let $1 \le u_0,u_1,u_\theta \le 2 $ and $1 \le v_0, v_1, 
v_\theta \le \infty$
 such that $1/{u_\theta} = (1-\theta)/{u_0} + \theta/{u_1}$ and 
$1/{v_\theta} = (1-\theta)/{v_0} + \theta/{v_1}$. Then 
$$
\lambda(\Pi_B,u_\theta,v_\theta) \le (1-\theta) \cdot \lambda(\Pi_B, 
u_0,v_0) + \theta \cdot \lambda(\Pi_B,u_1,v_1).
$$
\end{corollary}
\section{The proof of Theorem~\ref{khinchar}}
As a generalization of the notion of $\Lambda(p)$-sets, an orthonormal 
system $B \subset L_2(\mu)$ is said to be a {\em $\Lambda(p)$-system} 
if $B \subset L_p(\mu)$ and there exists a constant $c>0$ such that for 
all $f \in \Span B$ we have $\norm{f}_{L_p(\mu)} \le c \cdot 
\norm{f}_{L_2(\mu)}$; the infimum over all such constants $c$ is denoted
 by $K_p(B)$. Now one direction of the equivalence in Theorem~\ref{khinchar}
 can be formulated for general orthonormal systems:  
\begin{proposition}
\label{khintheo}
Let $B \subset L_2(\mu)$ be a 
$\Lambda(p)$-system for all $2 < p < \infty$. Then for all
 \mbox{$1 \le u,v \le \infty$} 
$$
 \lambda(\Pi_B,u,v) = \lambda(\Pi_\gamma,u,v) =
\begin{cases}
1/v &  \text{ if } 2 \le u \le \infty, \\
\max(0,1/2 +1/v-1/u) & \text{ if } 1 \le u \le 2.
\end{cases}
$$
\end{proposition}
\proof
 Although the limit order of $\Pi_\gamma$ is already known by the 
results of \cite{linde}, the following proof may also be used to compute
 it independently (at least the upper estimates; the lower ones are 
somehow simple), but for simplicity we fall back upon this knowledge. 
\par 
Since $\Pi_2 \subset \Pi_B \subset \Pi_\gamma$ (see \cite[4.15]{pw}; 
$\Pi_2$ denotes the Banach operator ideal of all $2$-summing operators), 
we only have to show that $\lambda(\Pi_B,u,v) \le \lambda(\Pi_\gamma,u,v)$,  
 and moreover, we conclude that 
$\lambda(\Pi_B,u,v)\le \lambda(\Pi_2,u,v) =\lambda(\Pi_\gamma,u,v)$ for all
 \mbox{$1 \le u \le \infty$} and $1 \le v \le 2$.
 For  $2 < v < \infty$ and $2 \le u \le \infty$ it can be easily seen 
 that 
$$
\Pi_B(\ell_u^m \Id \ell_v^m) \le K_v(B) \cdot m^{1/v}
$$
(just copy the proof of \cite[3.11.11]{pw}), hence---together with the continuity
 of the limit order, see \cite[14.4.8]{pietsch}---we obtain
$$
\lambda(\Pi_B,u,v) \le 
 1/v = \lambda(\Pi_\gamma,u,v)
$$ 
for all $2\le u,v\le \infty$. Now  the case 
$1 \le u \le 2 \le v \le \infty$ follows from  Corollary~\ref{interpollimit}:  
 For $1 \le u \le 2$ choose $u_0:=1$, $u_1:=2$, $v_0:=2$, $v_1:=\infty$, 
$\theta:=2/{u'}$ and 
$v_u$ such that  $1/{v_u} =1/u-1/2$. Then
$$
\lambda(\Pi_B,u,v_u) \le (1-\theta) \cdot \lambda(\Pi_B,1,2) 
+ \theta \cdot \lambda(\Pi_B,2,\infty) = 0.
$$
For arbitrary $1 \le u \le 2 \le v \le \infty$ 
factorize through $\ell_{v_u}^m$:
$$
\pi_B(\ell_u^m \Id \ell_v^m) \le m^{\max(0,1/v +1/2 -1/u)} \cdot 
\pi_B(\ell_u^m \Id \ell_{v_u}^m),
$$
hence $\lambda(\Pi_B,u,v) \le \max(0,1/v +1/2-1/u)= \lambda(\Pi_\gamma, 
u,v)$. 
\qed
\\[10pt]
The reverse implication in Theorem~\ref{khinchar} 
follows from \cite{baurdiss}: (b) trivially implies $\lambda(\Pi_B, \infty,
 \infty)=0$, and the comments after \cite[7.12]{baurdiss} then tell us that 
$\Pi_p \subset \Pi_\Lambda$ for all $2<p<\infty$. This in turn gives by 
\cite[9.6]{baurdiss} (see also \cite[5.1]{baur99}) 
that $\Lambda$ is a $\Lambda(p)$-set for all $2<p<\infty$. 
\par
Note that the last argument requires the setting of characters on a 
compact abelian group; Baur has recently informed us that her results are
 also valid for the non-abelian case, and therefore our 
Theorem~\ref{khinchar} as well. 
\section{The proof of Theorem~\ref{gammaschatten}}
\proof
For $1 \le v \le 2$ by \cite{tom} $\ui_v$ is of cotype~$2$, hence  
$$
 n^{1/v + \min(1/2,1-1/u)}  = \pi_2(\ui_u^n \Id \ui_v^n) 
\ge \pi_\gamma(\ui_u^n \Id \ui_v^n) 
\ge \cotype{2}(\ui_v)^{-1} \cdot  \pi_2(\ui_u^n \Id \ui_v^n) 
$$   
(see e.\,g. \cite[Corollary 3]{dm98}) for $1 \le u \le \infty$ and 
$1 \le v \le 2$; here $\cotype{2}(\ui_v)$ denotes the Gaussian 
cotype 2 constant of $\ui_v$. We are left with 
 the case $2 \le v \le \infty$; first let $u=v=\infty$. 
Then by \cite[4.15.18]{pw} 
(a result of Pisier, see also \cite[4.4]{pisier89} and \cite{pisierLN}) 
and \cite[3.3]{figiel} for each $\varepsilon>0$
$$
\pi_\gamma(\ui_\infty^n \Id \ui_\infty^n) \asymp \sqrt{D(\ui_\infty^n,
\varepsilon)} 
\asymp n^{1/2},
$$
where $D(X,\varepsilon)$ denotes the Dvoretzky
 dimension of a Banach space $X$, i.\,e. the largest $m$ such that there 
exists an $m$-dimensional subspace $X_m$ of $X$ with Banach--Mazur 
distance $d(X_m,\ell_2^m) \le
 1+\varepsilon$ (see \cite[4.15.15]{pw}). Now the general case 
$2 \le u,v \le \infty$ follows by factorization:
$$
\pi_\gamma(\ui_u^n \Id \ui_v^n) \le 
n^{1/v} \cdot \pi_\gamma(\ui_\infty^n \Id \ui_\infty^n) \asymp n^{1/v+1/2},
$$
and conversely
$$
\pi_\gamma(\ui_u^n \Id \ui_v^n) \ge n^{1/v-1/2} \cdot 
\pi_\gamma(\ui_2^n \Id \ui_2^n) = n^{1/v+1/2}.
$$
The case $1 \le u \le 2 \le v \le \infty$ is done by interpolation: 
We have (recall that \mbox{$\pi_2(\ui_1^n \Id \ui_2^n)=n^{1/2}$})
$$
\pi_\gamma(\ui_1^n \Id \ui_2^n) \asymp \pi_\gamma(\ui_2^n \Id \ui_\infty^n)
 \asymp n^{1/2},
$$
hence for $1 < u <2 < v_u < \infty$ and $0 < \theta <1$ such that 
$1/{v_u} =1/u -1/2$ and $\theta=2/{u'}$
$$
\pi_\gamma(\ui_u^n \Id \ui_{v_u}^n) \le d_\theta[\ui_1^n,\ui_2^n] \cdot 
\pi_\gamma(\ui_1^n \Id \ui_2^n)^{1-\theta} \cdot \pi_\gamma(\ui_2^n \Id 
\ui_\infty^n)^\theta \asymp n^{1/2};
$$
recall that $\sup_n d_\theta[\ui_1^n,\ui_2^n] < \infty$ by \eqref{dtheta2}, 
 and that $[\ui_1^n,\ui_2^n]_\theta = \ui_u^n$ and 
$[\ui_2^n, \ui_\infty^n]_\theta = \ui_{v_u}^n$ hold isometrically (this 
can be deduced from e.\,g. \cite[Satz 8]{pt} and the complex reiteration
 theorem \cite[4.6.1]{BL}).
 The remaining estimates now follow easily from
$$
\pi_\gamma(\ui_u^n \Id \ui_v^n) \ge \pi_\gamma(\ell_2^n \Id \ell_2^n)=
n^{1/2},
$$ 
$$
\pi_\gamma(\ui_u^n \Id \ui_v^n) \ge \pi_\gamma(\ui_u^n \Id \ui_2^n) 
\cdot \norm{\ui_v^n \Id \ui_2^n}
$$
and
$$
\pi_\gamma(\ui_u^n \Id \ui_v^n) \le \pi_\gamma(\ui_u^n \Id \ui_{v_u}^n) 
\cdot \norm{\ui_{v_u}^n \Id \ui_v^n}.
$$
\qed
\par
The contents of this article are part of the author's thesis written
 at the Carl von Ossietzky University of Oldenburg under the supervision
 of Prof. Dr. Andreas Defant.
\providecommand{\bysame}{\leavevmode\hbox to3em{\hrulefill}\thinspace}

\end{document}